\documentclass[11pt]{amsart}

\usepackage[T1]{fontenc}
\usepackage[utf8]{inputenc}
\usepackage{amsmath,amssymb,amsthm}
\usepackage{booktabs,longtable,array}
\usepackage[margin=2.5cm]{geometry}
\usepackage{graphicx}

\newtheorem{theorem}{Theorem}
\newtheorem{lemma}{Lemma}
\newtheorem{corollary}{Corollary}
\theoremstyle{remark}
\newtheorem*{remark}{Remark}

\newcommand{\Sol}{\v{S}olt\'es}

\title[An infinite family with many \v{S}olt\'es vertices]{An infinite family of graphs with \v{S}olt\'es fraction\\ uniformly above one half}

\author{A.~A.~Dobrynin}
\address{Sobolev Institute of Mathematics, Siberian Branch of the Russian Academy of Sciences, Novosibirsk, 630090, Russia}
\email{dobr@math.nsc.ru}


\begin{document}

\begin{abstract}
The Wiener index $W(G)$ of a connected graph $G$ is the sum of the distances between all unordered pairs of vertices.
A vertex $x$ is a \Sol{} vertex if $G-x$ is connected and $W(G-x)=W(G)$.
We construct an infinite family of graphs in which the proportion of \Sol{} vertices is strictly greater than
$(\sqrt{21}-3)/3=0.5275\ldots>1/2$ and decreases to this value from $2/3$.
The family is also described by a generalized Pell equation and starts from a $12$-vertex example of Akhmejanova et al. (Discrete Appl. Math. 325 (2023) 37--42).
We also show that the construction remains valid after adding edges among certain vertices.
\end{abstract}

\maketitle

\section*{Introduction}

In this paper we are concerned with finite simple undirected
connected graphs. The vertex set of $G$ is
denoted by $V(G)$. 
The cardinality of the set $V(G)$ is called the
\emph{order} of the graph $G$.
If $u$ and $v$ are vertices of $G$, then the number of
edges in a shortest path connecting them is said to be
their distance and is denoted by $d(u,v)$.
The \emph{transmission} of a vertex $v$ in a graph $G$ is the sum of distances between $v$ and all
other 
vertices of $G$, $T_G(v)=\sum_{u\in V(G)}d_G(v,u)$. 
The Wiener index of $G$ is a distance-based graph invariant defined as the sum of distances 
between all unordered pairs of vertices of $G$, 
$$
W(G)=\sum_{\{u,v\}\subseteq V(G)}d_G(u,v)
=\frac12\sum_{v\in V(G)}T_G(v).
$$

It was originally introduced as a structural descriptor for alkanes,
whose carbon skeletons are trees with vertex degrees at most four \cite{Wiener}.
Mathematical properties and chemical applications of the
Wiener index have been intensively studied over the past
decades 
\cite{BonchevRouvray1991,DobryninEntringerGutman, DobryninGutmanKlavzarZigert, Gutman1986, KnorSkrekovskiTepeh2016, KnorSkrekovskiTepeh2024,Trinajstic1992,WagnerWang2018}.

One of the research directions concerning the index involves 
describing its properties under graph operations. A vertex $v\in V(G)$ is called a \emph{\v{S}olt\'{e}s vertex} if $G-v$ is connected and $W(G-v)=W(G)$.
In 1991, \v{S}olt\'es~\cite{Soltes} posed the problem of finding all graphs in which every vertex is a \Sol{} vertex.
The cycle $C_{11}$ is the only known example, and the problem remains open.

Since the full problem seems out of reach, various relaxations have been studied.
Subsequent work has investigated preservation upon the removal of particular vertices, 
including constructions in the classes of unicyclic, cactus, cubic, and regular graphs 
\cite{BasicKnorSkrekovski, BokJedlickovaMaxova2019, BokJedlickovaMaxova2021, Dobrynin2023, KnorMajstorovicSkrekovskiAMC, KnorMajstorovicSkrekovskiDAM}.
The problem of finding graphs with arbitrarily many \Sol{} vertices was posed in \cite{KnorMajstorovicSkrekovskiAMC}.
For every $k\ge1$, infinitely many graphs with at least $k$ \Sol{} vertices were constructed in \cite{BokJedlickovaMaxova2019,BokJedlickovaMaxova2021,HuZhuWuShaoFahad2021}.
Such families of cubic 2-connected graphs were constructed in \cite{BasicKnorSkrekovski}.
Preservation of the Wiener index under the simultaneous deletion of $k\ge2$ vertices was studied in \cite{HuZhuWuShaoFahad2021}, following a problem posed in \cite{KnorMajstorovicSkrekovskiAMC}.

Other relaxations concern graphs in which at least a fixed proportion of the vertices are \Sol{} vertices~\cite{BasicKnorSkrekovski}.
We measure this proportion by the \emph{\Sol{} fraction} $R_0(G)=N_S/|V(G)|$,
where $N_S$ is the number of \Sol{} vertices.
Two values of $R_0$ are of particular interest. The first is $1/2$.
It is difficult to find graphs whose \Sol{} fraction is close to or greater than this value.
The constructions in \cite{DobryninVorobev} produce three graphs with $R_0\in\{4/7,5/9,6/11\}$.
A graph of order $12$ with \Sol{} fraction exactly $1/2$ is shown in \cite{Cambie2026}.
It can be verified computationally that some generalized Petersen graphs $GP(n,k)$ have \Sol{} fraction exactly $1/2$.
Infinite families of graphs whose \Sol{} fraction tends to $1/2$ from below
were constructed in \cite{AkhmejanovaEtAl,DobryninVorobev}.
The second value is $2/3$. Only a few graphs with this \Sol{} fraction are known.
A graph of order $12$ with \Sol{} fraction $2/3$ appears in~\cite{AkhmejanovaEtAl}.
Eight graphs with \Sol{} fraction $2/3$, including this one, are described by Cambie~\cite{Cambie2026}. 
To the best of our knowledge, no infinite family of graphs with \Sol{} fraction above $1/2$ has previously been presented.

In this note we construct an explicit infinite family of graphs whose \Sol{} fractions have a common lower bound strictly greater than $1/2$.
The first member is the graph of order $12$ with \Sol{} fraction $2/3$ from \cite{AkhmejanovaEtAl}.
Further infinite families are obtained by adding edges between certain vertices without changing the set of \Sol{} vertices.

\section{Main result}\label{sec:construction}

\begin{figure}[ht]
\centering
\includegraphics[width=\linewidth]{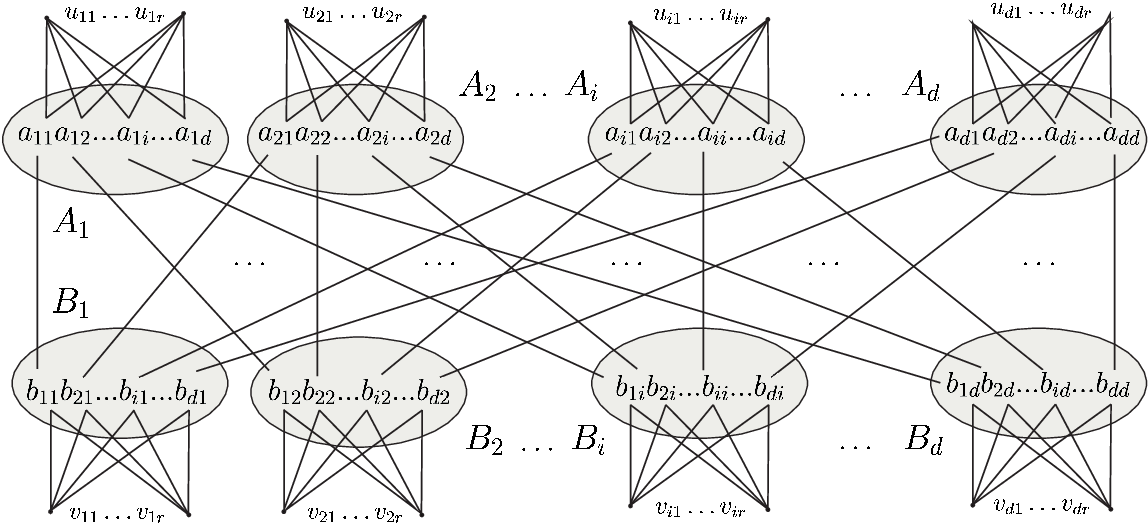}
\caption{The graph $G(d,r)$.}\label{fig:construction}
\end{figure}

For integers $d\ge 2$ and $r\ge 1$, consider the graph $G(d,r)$ shown in Fig.~\ref{fig:construction}.
Its vertex set consists of \emph{base vertices} $a_{ij}$, $b_{ij}$ and \emph{auxiliary vertices} $u_{it}$, $v_{jt}$,
where $1\le i,j\le d$ and $1\le t\le r$.
For each $i$, the vertices $a_{i1},\dots,a_{id}$ induce a clique $A_i$ of order $d$.
For each $j$, the vertices $b_{1j},\dots,b_{dj}$ induce a clique $B_j$ of order $d$.
The cliques $A_i$ and $B_j$ are called base cliques.
These $2d$ cliques are pairwise disjoint and are joined by the edges $a_{ij}b_{ij}$.
Each auxiliary vertex $u_{it}$ is adjacent to all vertices of $A_i$, and each $v_{jt}$ to all vertices of $B_j$.

Thus $G(d,r)$ has $2d^2$ base vertices of degree $d+r$, and $2dr$ auxiliary vertices of degree $d$. 
The order of $G(d,r)$ is $2d(d+r)$.
Deleting the auxiliary vertices and contracting each of the cliques $A_i$, $B_j$ to a single vertex gives the complete bipartite graph $K_{d,d}$.
The base vertices and the auxiliary vertices form the two orbits of the automorphism group of $G(d,r)$.
Indeed, this group contains the maps that permute the first indices $i$ (together with the corresponding $u_{it}$),
those that permute the second indices $j$ (together with the corresponding $v_{jt}$),
those that permute the index $t$,
and the involution $a_{ij}\leftrightarrow b_{ji}$, $u_{it}\leftrightarrow v_{it}$.

Using $G(d,r)$, we construct an infinite family of graphs in which more than half of the vertices are \Sol{} vertices.

\newpage
\begin{theorem}\label{thm:main}
Let $G_n=G(d_n,r_n)$, where $d_n=2s_n$, $s_0=r_0=1$, and
\[
s_{n+1}=14s_n+5r_n-4,\qquad r_{n+1}=25s_n+9r_n-8\qquad(n\ge 0).
\]
Then $G_n$ of order $2d_n(d_n+r_n)$  has exactly $2d_n^2$ \Sol{} vertices, 
\[
R_0(G_n)=\frac{d_n}{d_n+r_n}>\frac{\sqrt{21}-3}{3}=0.527525\ldots>\frac12 ,
\]
and $R_0(G_n)$ decreases strictly from $2/3$ to $(\sqrt{21}-3)/3$.
\end{theorem}

The construction of $G_n$ allows us to obtain new families of graphs with a similar structure satisfying Theorem~\ref{thm:main}.

\begin{corollary}\label{cor:main}
Let $G_n'$ be obtained from $G_n$ by adding arbitrary edges within each set of $r_n$ auxiliary vertices
attached to the same base clique.
Then $G_n'$ has the same order and the same \Sol{} vertices as $G_n$, so Theorem~\ref{thm:main} holds for $G_n'$.
\end{corollary}

\section{Changes in distances after vertex removal}\label{sec:distance}

Let $G=G(d,r)$ and $\Delta_G(x)=W(G)-W(G-x)$ for $x\in V(G)$.

\begin{lemma}\label{lem:base}
If $x$ is a base vertex of $G$, then
\begin{equation}\label{eq:delta}
\Delta_G(x)=5d^2-8d+4-2dr-4r^2 .
\end{equation}
\end{lemma}

\begin{proof}
The base vertices form a single orbit, so we may take $x=a_{11}$.
A shortest path between base vertices contains no auxiliary vertex, since the neighbours of every auxiliary vertex form a clique.

The vertex $x$ has $d$ base neighbours: $b_{11}$ and $a_{1j}$ for $2\le j\le d$.
At distance $2$ from $x$ there are $2(d-1)$ base vertices: $b_{1j}$ for $2\le j\le d$ and $b_{i1}$ for $2\le i\le d$.
At distance $3$ from $x$ there are $d(d-1)$ base vertices, forming two sets:
\begin{itemize}
\item the $d-1$ vertices $a_{i1}$, $2\le i\le d$, reached by the path $(a_{11},b_{11},b_{i1},a_{i1})$;
\item the $(d-1)^2$ vertices $b_{ij}$, $2\le i,j\le d$, reached by the path $(a_{11},a_{1j},b_{1j},b_{ij})$.
\end{itemize}
Neither set meets the first two distance layers, so these paths are shortest.
At distance $4$ from $x$ there are $(d-1)^2$ base vertices $a_{ij}$, $2\le i,j\le d$, reached by the path
$(a_{11},a_{1j},b_{1j},b_{ij},a_{ij})$.
This path is shortest because $a_{ij}$ has no base neighbour at distance less than $3$ from $x$.

For each $t\in\{1,\dots,r\}$, the auxiliary vertices $u_{1t}$, $v_{1t}$, $v_{jt}$ $(2\le j\le d)$ and $u_{it}$ $(2\le i\le d)$
lie at distances $1,2,3,4$ from $x$, respectively.
Hence
\begin{equation}\label{eq:Tx}
T_G(x)=7d^2-6d+(7d-4)r .
\end{equation}

Table~\ref{tab:paths} lists all unordered pairs of vertices of $G-x$ whose distance increases when $x=a_{11}$ is deleted.
Each path is written as a sequence of vertices; $i,j\in\{2,\dots,d\}$ and $s,t\in\{1,\dots,r\}$ vary independently.
The index $2$ appearing in some paths is admissible because $d\ge 2$.

\begin{table}[ht]
\centering
\caption{Distances and shortest paths before and after deleting $a_{11}$.}\label{tab:paths}
\small
\setlength{\tabcolsep}{4pt}
\renewcommand{\arraystretch}{1.25}
\begin{tabular}{@{}l c l c l c@{}}
\toprule
Pair & Before & Path in $G$ & After & Path in $G-x$ & Pairs\\
\midrule
$a_{1j},b_{i1}$ & 3 & $(a_{1j},a_{11},b_{11},b_{i1})$ & 5 & $(a_{1j},b_{1j},b_{ij},a_{ij},a_{i1},b_{i1})$ & $(d-1)^2$\\
$a_{1j},b_{11}$ & 2 & $(a_{1j},a_{11},b_{11})$ & 6 & $(a_{1j},b_{1j},b_{2j},a_{2j},a_{21},b_{21},b_{11})$ & $d-1$\\
$b_{11},b_{1j}$ & 3 & $(b_{11},a_{11},a_{1j},b_{1j})$ & 5 & $(b_{11},b_{21},a_{21},a_{2j},b_{2j},b_{1j})$ & $d-1$\\ \hline
$a_{i1},u_{1t}$ & 4 & $(a_{i1},b_{i1},b_{11},a_{11},u_{1t})$ & 5 & $(a_{i1},a_{i2},b_{i2},b_{12},a_{12},u_{1t})$ & $(d-1)r$\\
$a_{1j},v_{1t}$ & 3 & $(a_{1j},a_{11},b_{11},v_{1t})$ & 6 & $(a_{1j},b_{1j},b_{2j},a_{2j},a_{21},b_{21},v_{1t})$ & $(d-1)r$\\
$b_{i1},u_{1t}$ & 3 & $(b_{i1},b_{11},a_{11},u_{1t})$ & 6 & $(b_{i1},a_{i1},a_{i2},b_{i2},b_{12},a_{12},u_{1t})$ & $(d-1)r$\\
$b_{1j},v_{1t}$ & 4 & $(b_{1j},a_{1j},a_{11},b_{11},v_{1t})$ & 5 & $(b_{1j},b_{2j},a_{2j},a_{21},b_{21},v_{1t})$ & $(d-1)r$\\
$b_{11},u_{1t}$ & 2 & $(b_{11},a_{11},u_{1t})$ & 7 & $(b_{11},b_{21},a_{21},a_{22},b_{22},b_{12},a_{12},u_{1t})$ & $r$\\
$b_{11},v_{jt}$ & 4 & $(b_{11},a_{11},a_{1j},b_{1j},v_{jt})$ & 5 & $(b_{11},b_{21},a_{21},a_{2j},b_{2j},v_{jt})$ & $(d-1)r$\\ \hline
$u_{1t},v_{1s}$ & 3 & $(u_{1t},a_{11},b_{11},v_{1s})$ & 7 & $(u_{1t},a_{12},b_{12},b_{22},a_{22},a_{21},b_{21},v_{1s})$ & $r^2$\\
\bottomrule
\end{tabular}
\end{table}

To see that the listed paths are shortest and that the list is exhaustive, consider the subgraph induced by the base vertices and contract each base clique to a single vertex.
Deleting $a_{11}$ removes the edge $a_{11}b_{11}$ joining the cliques $A_1$ and $B_1$ (see Fig.~\ref{fig:construction}).
In $K_{d,d}$, deleting this edge changes only the distance between its endpoints, from $1$ to $3$.
In the subgraph induced by the base vertices, consecutive edges joining distinct cliques are separated by an edge inside a clique.
In $G-a_{11}$, every base neighbour of $b_{11}$ lies in the clique $B_1$ containing $b_{11}$.
Together, these facts justify the first three rows of Table~\ref{tab:paths}.

\newpage

For a base vertex $y$ and an auxiliary vertex $z$ of $H\in\{G,G-x\}$,
\[
d_H(y,z)=1+\min_{a\in N_H(z)}d_H(y,a),
\]
where $N_H(z)$ denotes the set of neighbours of $z$ in $H$.
The distances in the next six rows follow from this formula.
For two distinct auxiliary vertices $z$ and $z'$,
\[
d_H(z,z')=2+\min_{\substack{a\in N_H(z)\\ b\in N_H(z')}}d_H(a,b).
\]
The last row follows from this identity.
Distances within a base clique and between auxiliary vertices attached to the same base clique do not change. The total increase of distances among the remaining pairs is therefore
\[
D(x)=\sum_{\{y,z\}\subseteq V(G)\setminus\{x\}}\bigl(d_{G-x}(y,z)-d_G(y,z)\bigr)
=2(d-1)(d+2)+(9d-4)r+4r^2 .
\]
Since $\Delta_G(x)=T_G(x)-D(x)$, equation~\eqref{eq:Tx} yields~\eqref{eq:delta}.
\end{proof}

If $x$ is an auxiliary vertex, its neighbours form a clique, so deleting $x$ does not change any distance between the remaining vertices.
Hence $\Delta_G(x)=T_G(x)>0$.
For an auxiliary vertex $z$,
\[
T_G(z)=7d^2-4d+(8d-3)r-2 .
\]
Since the Wiener index is half the sum of the vertex transmissions, we obtain
\[
W(G(d,r))=7d^4-6d^3+(14d^3-8d^2-2d)\,r+d(8d-3)\,r^2 .
\]

\section{Proof of Theorem~\ref{thm:main}}\label{sec:recurrence}

Let $d=2s$.
By equation~\eqref{eq:delta}, $\Delta_G(x)=4F(s,r)$ for a base vertex $x$, where
\[
F(s,r)=5s^2-4s+1-sr-r^2 .
\]
Direct expansion gives the polynomial identity
\begin{equation}\label{eq:identity}
F(14s+5r-4,\;25s+9r-8)=F(s,r).
\end{equation}
Since $F(s_0,r_0)=F(1,1)=0$, identity~\eqref{eq:identity} and induction give $F(s_n,r_n)=0$ for all $n\ge0$.
The recurrences yield $s_{n+1}>s_n\ge1$ and $r_{n+1}>r_n\ge1$, so the orders of the graphs $G_n$ strictly increase.
By Lemma~\ref{lem:base}, every base vertex $x$ of $G_n$ satisfies
$\Delta_{G_n}(x)=4F(s_n,r_n)=0$, whereas $\Delta_{G_n}(z)>0$ for every auxiliary vertex $z$.
Consequently, the \Sol{} vertices of $G_n$ are precisely its $2d_n^2$ base vertices.

The equation $F(s,r)=0$, regarded as a quadratic equation in $r$, has the positive root
\[
r=\frac{-s+\sqrt{21s^2-16s+4}}{2}.
\]
Therefore
\begin{equation}\label{eq:R0}
R_0(G_n)=\frac{2d_n^2}{2d_n(d_n+r_n)}=\frac{2s_n}{2s_n+r_n}=\frac{4}{3+\sqrt{21-16/s_n+4/s_n^2}} .
\end{equation}
For $s\ge1$, the expression $21-16/s+4/s^2$ is strictly increasing in $s$ and less than $21$.
Since $s_n\to\infty$, the fractions in~\eqref{eq:R0} decrease strictly to
\[
\frac{4}{3+\sqrt{21}}=\frac{\sqrt{21}-3}{3}=0.527525\ldots>\frac12 .
\]
This completes the proof of Theorem~\ref{thm:main}.

\begin{figure}[ht]
\centering
\includegraphics[width=0.4\linewidth]{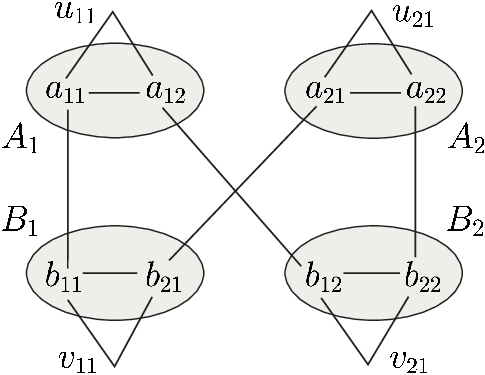}
\caption{The graph $G_0=G(2,1)$.}\label{fig:example}
\end{figure}
 
The first graph $G_0=G(2,1)$ in the family is shown in Fig.~\ref{fig:example}.
It has order $12$ and $R_0(G(2,1))=2/3$. 
This graph was first introduced by Akhmejanova et al.~\cite{AkhmejanovaEtAl}.

\begin{remark}[A Pell-type interpretation]
By Lemma~\ref{lem:base}, deleting a base vertex of $G(2s,r)$ preserves the Wiener index if and only if $F(s,r)=0$.
Thus, to construct infinitely many such graphs, it suffices to find infinitely many positive integer solutions $(s,r)$ of this equation.
An alternative way to obtain these solutions is to use a generalized Pell equation (see \cite{Barbeau2003} for background).
The substitution $x=42s-16$, $y=2s+4r$ transforms $F(s,r)=0$ into
\[
x^2-21y^2=-80 .
\]
The initial pair $(s_0,r_0)=(1,1)$ corresponds to the solution $(x_0,y_0)=(26,6)$.
Starting with this solution, we multiply $x_n+y_n\sqrt{21}$ by $(23+5\sqrt{21})/2$ and define the next pair by
\begin{equation}\label{eq:pell-step}
x_{n+1}+y_{n+1}\sqrt{21}=\frac{23+5\sqrt{21}}{2}\,\bigl(x_n+y_n\sqrt{21}\bigr).
\end{equation}
The new pair satisfies the same Pell equation.
Indeed, multiplying~\eqref{eq:pell-step} by its conjugate, we obtain
\[
x_{n+1}^{2}-21y_{n+1}^{2}
=\frac{23^{2}-21\cdot5^{2}}{4}\bigl(x_n^{2}-21y_n^{2}\bigr)
=x_n^{2}-21y_n^{2}=-80.
\]
Comparing the coefficients of $1$ and $\sqrt{21}$ in~\eqref{eq:pell-step}, we obtain
\[
x_{n+1}=\frac{23x_n+105y_n}{2},\quad y_{n+1}=\frac{5x_n+23y_n}{2}.
\]
Returning to the original parameters, we obtain
\[
\begin{aligned}
s_{n+1}&=\frac{x_{n+1}+16}{42}=14s_n+5r_n-4,\\
r_{n+1}&=\frac{y_{n+1}-2s_{n+1}}{4}=25s_n+9r_n-8.
\end{aligned}
\]
These are the recurrences in Theorem~\ref{thm:main}.
They ensure that $s_n$ and $r_n$ are positive integers and that $s_{n+1}>s_n$.
Hence they generate infinitely many distinct graphs $G_n=G(2s_n,r_n)$.
\end{remark}

Table~\ref{tab:seq} contains parameters of the smallest graphs $G_n=G(d_n, r_n)$.

\begin{table}[ht]
\centering
\caption{Initial members of the infinite sequence.}\label{tab:seq}
\small
\renewcommand{\arraystretch}{1.2}
\begin{tabular}{@{}r r r r r l@{}}
\toprule
$n$ & $d_n$ & $r_n$ & Order & \Sol{} vertices & $R_0$\\
\midrule
0 & 2 & 1 & 12 & 8 & $2/3$\\
1 & 30 & 26 & 3360 & 1800 & $15/28\approx 0.535714$\\
2 & 672 & 601 & 1710912 & 903168 & $672/1273\approx 0.527887$\\
3 & 15410 & 13801 & 900283020 & 474936200 & $15410/29211\approx 0.527541$\\
4 & 353742 & 316826 & 474416130912 & 250266805128 & $176871/335284\approx 0.527526$\\
5 & 8120640 & 7273201 & 250015681956480 & 131889588019200 & $8120640/15393841\approx 0.527525$\\
\bottomrule
\end{tabular}
\end{table}

\section{Proof of Corollary~\ref{cor:main}}\label{sec:corollary}

We first establish the following lemma.

\begin{lemma}\label{lem:twins}
Let $G$ be a connected graph, and let $Z$ be an independent set of vertices that all have the same neighbourhood $X$, where $|X|\ge2$.
Let $H$ be obtained from $G$ by adding arbitrary edges between vertices of $Z$.
Then, for every $x\in V(G)\setminus Z$ such that $G-x$ is connected,
\[
W(H)-W(H-x)=W(G)-W(G-x).
\]
\end{lemma}

\begin{proof}
Let $E^+=E(H)\setminus E(G)$ and $m=|E^+|$.
Obviously, a shortest path in $H$ with at least one endpoint outside $Z$ contains no added edge.
Since any two distinct vertices of $Z$ are at distance $2$ in $G$, for distinct $u,v\in V(G)$ we have $d_G(u,v)-d_H(u,v)=1$ if $uv\in E^+$, and $0$ otherwise.
Hence $W(H)=W(G)-m$. Since $x\notin Z$, all $m$ added edges remain in $H-x$.
Since $|X|\ge2$, the vertices of $Z$ always have a common neighbour in $G-x$, and the same argument yields $W(H-x)=W(G-x)-m$.
Subtracting the two equalities proves the lemma.
\end{proof}

\begin{proof}[Proof of Corollary~\ref{cor:main}]
Apply Lemma~\ref{lem:twins} successively to the $2d_n$ sets of $r_n$ auxiliary vertices attached to the same base clique.
For each of them, $X$ is that clique and $|X|=d_n\ge2$.
For every base vertex $x$ this implies
\[
\Delta_{G_n'}(x)=\Delta_{G_n}(x)=0 .
\]
If $z$ is an auxiliary vertex, deleting $z$ does not change distances between the remaining vertices.
Hence $\Delta_{G_n'}(z)>0$.
Consequently, $G_n'$ has the same \Sol{} vertices as $G_n$.
\end{proof}

As a particular case, add all edges within each set of $r$ auxiliary vertices, and denote the resulting graph by $G^{c}(d,r)$.
It has the same \Sol{} vertices as $G(d,r)$, and
\[
W\bigl(G^{c}(d,r)\bigr)=W\bigl(G(d,r)\bigr)-dr(r-1).
\]
Each base clique, together with its auxiliary vertices, now induces $K_{d+r}$.
In $G^{c}(d,r)$, base vertices have degree $d+r$, and auxiliary vertices have degree $d+r-1$.

\begin{table}[ht]
\centering
\caption{Two versions of the graph of order $3360$.}\label{tab:versions}
\small
\renewcommand{\arraystretch}{1.2}
\begin{tabular}{@{}l r r r r@{}}
\toprule
Auxiliary vertices & Base vertex degree & Auxiliary vertex degree & $W$ & \Sol{} vertices\\
\midrule
Independent & 56 & 30 & 19953600 & 1800\\
Cliques & 56 & 55 & 19934100 & 1800\\
\bottomrule
\end{tabular}
\end{table}

\section*{Conclusion}
We have constructed an infinite family of graphs in which the proportion of \Sol{} vertices is strictly greater than
$(\sqrt{21}-3)/3=0.5275\ldots>1/2$ and decreases to this value from $2/3$.
Adding arbitrary edges within each set of auxiliary vertices attached to the same base clique preserves the set of \Sol{} vertices.
The distance formulas for $G(d,r)$ were verified computationally for several small values of $d$ and $r$.

\section*{Acknowledgements}
This work was supported by the state contract of the Sobolev Institute of Mathematics (project number FWNF-2022-0017).

\section*{Declaration of generative AI use}
ChatGPT (OpenAI) was used to assist with mathematical discussions and computational checks.
The author reviewed and revised the material and takes full responsibility for the final content.

\end{document}